\documentclass[reqno,11pt]{amsart}
\usepackage{url}
\usepackage{latexsym}
\usepackage{lscape} 
\usepackage{geometry}
\usepackage{pdflscape}
\usepackage{float}
\usepackage{amssymb}
\usepackage{amscd}
\usepackage{graphicx,color,cancel}
\usepackage{tikz, tikz-cd}
\usepackage[all]{xy}
\usepackage{subfiles} 
\usepackage[extrasp=0em,mono]{inconsolata}  
\usepackage{ytableau}

\usepackage{amssymb,amscd,verbatim, xcolor, amsthm,amsmath,amsgen,xspace,bm}
\usepackage{hyperref}
\usepackage{soul,todonotes}
\usepackage[normalem]{ulem}
\usepackage[mathscr]{eucal}
\usepackage{dutchcal}
\usepackage{graphicx,color,cancel}
\usepackage{xcolor}

\newtheorem{thm}[equation]{Theorem}

\newtheorem{ex}[equation]{Example}
\newtheorem{lem}[equation]{Lemma}
\newtheorem{cor}[equation]{Corollary}

\newtheorem{prop}[equation]{Proposition}

\theoremstyle{remark}
\newtheorem{rem}[equation]{Remark}

\theoremstyle{definition}

\numberwithin{equation}{section}

\newcommand{\mult}{\operatorname{mult}}

\newcommand{\cal}{\mathcal}

\newcommand{\bbar}{\,\big|\,}
\newcommand{\sbar}{\,|\,}

\newcommand\eps{{\varepsilon}}
\newcommand\la{{\lambda}}

\newcommand\ot{\otimes}

\newcommand\GL{\operatorname{GL}}

\newcommand\ini{\operatorname{in}}
\newcommand\Spec{\operatorname{Spec}}
\newcommand\LCM{\operatorname{LCM}}
\newcommand\WDF{\operatorname{WDF}}

\newcommand{\pf}{\begin{proof}}
\newcommand{\epf}{\end{proof}}
\newcommand{\eq}{\begin{equation}}
\newcommand{\eeq}{\end{equation}}
\newcommand{\eqn}{\begin{equation*}}
\newcommand{\eeqn}{\end{equation*}}

\newcommand{\frgl}{\mathfrak{gl}}

\newcommand{\bbC}{\mathbb{C}}

\newcommand{\bbN}{\mathbb{N}}
\newcommand{\bbQ}{\mathbb{Q}}

\newcommand{\bbZ}{\mathbb{Z}}

\newcommand{\Pol}{\mathcal{P}}

\newcommand{\cA}{\mathcal{A}}

\newcommand{\wt}{\operatorname{wt}}
\newcommand{\LT}{\operatorname{LT}}

\newcommand{\smat}{\left(\begin{smallmatrix}}
\newcommand{\esmat}{\end{smallmatrix}\right)}
\newcommand{\pmat}{\begin{pmatrix}}
\newcommand{\epmat}{\end{pmatrix}}
\newcommand{\mat}{\begin{matrix}}
\newcommand{\emat}{\end{matrix}}

\newcommand{\bbP}{\mathbb{P}}

\begin{document}

\title{Toric degeneration for hybrid qudit spaces}
\author{Jing-Song Huang}
\address[Huang]{School of Science and Engineering, The Chinese Universit of Hong Kong, Shenzhen, 2001 Long Xiang Road, Long Gang District, Shenzhen, Guangdong, China.}
\email{huangjingsong@cuhk.edu.cn}
\author{Soo Teck Lee}
\address[Lee]{Department of Mathematics\\
National University of Singapore\\
2 Science Drive 2\\
Singapore 117543, Singapore.} \email {matleest@nus.edu.sg}

\author{Pavle Pand\v zi\'c}
\address[Pand\v zi\'c]{Department of Mathematics, Faculty of Science, University of Zagreb, Bijeni\v cka 30, 10000 Zagreb, Croatia.}
\email{pandzic@math.hr}
\thanks{J.-S.~Huang is supported by NSFC Grants No. 12271460 and 12341101, by Guangdong Provincial Grants 2024A1515011456 and GDZX2403006, and by 	
	Shenzhen Municipal Grants JCYJ20241202124023031.  S.~T.~Lee is supported by NUS grant A-8002916-00-00. 
	P.~Pand\v zi\'c is supported by "Implementation of cutting-edge research and its application as part of the Scientific Center of Excellence for Quantum and Complex Systems, and Representations of Lie Algebras", Grant No. PK.1.1.10.0004, co-financed by the European Union through the European Regional Development Fund - Competitiveness and Cohesion Programme 2021- 2027., by the Croatian Science Foundation (HRZZ), grant no. IP-2025-02-6514, and
	by the European Union - NextGenerationEU through the National Recovery and Resilience Plan 2021-2026. Institutional grant of University of Zagreb Faculty of Science (IK IA 1.1.3. Impact4Math).}
\begin{abstract} We consider covariant algebras attached to certain hybrid qudit spaces. We describe these algebras in terms of generators and relations, and we show they are flat deformations of certain explicitly described semigroup algebras.
\end{abstract}

\maketitle

\section{Introduction}

The problem of understanding the decomposition under $GL(n_1)\times\dots\times GL(n_k)$ of the algebra of polynomials on the tensor product of standard modules $\bbC^{n_1}\otimes\dots\otimes\bbC^{n_k}$ is very old and well known. It is related to various problems in geometry, physics and combinatorics. For example, the Kronecker coefficients, i.e., multiplicities in the decomposition of tensor products of representations of the  symmetric group, can be recovered from the $GL(n)\times GL(n)\times GL(n)$-decomposition of polynomials on $\bbC^n\otimes\bbC^n\otimes\bbC^n$; see for example \cite{BV}. 

In general, this is a  very difficult problem, except for the case $k=2$, which is very well known, see e.g. \cite{H}. 
The interest in some other special cases has surged recently, driven by relations with the rapidly developing subject of quantum computing. Thus the study of $k$-qubit states is directly related to the above problem for $n_1=\dots =n_k=2$, while qudit states are related to the already mentioned case of $k=3$, $n_1=n_2=n_3=n$. Some references to mention here are \cite{La}, \cite{W}, \cite{BLT}, \cite{Bri}, \cite{Lu}, \cite{WDGC},  
\cite{B}, \cite{BB}, \cite{MW1}, \cite{MW2} and \cite{FH}.

To understand the decomposition of the algebra $\cal P=\cal P(\bbC^{n_1}\otimes\dots\otimes\bbC^{n_k})$ of polynomials on $\bbC^{n_1}\otimes\dots\otimes\bbC^{n_k}$  under $GL(n_1)\times\dots\times GL(n_k)$, one would like to describe the algebra of covariants, i.e., the algebra of highest weight vectors of the irreducible submodules of $\cal P$. For the case of 3 qubits, such a description by generators and relations was described as early as  1881 in \cite{LP}. In larger cases the analysis rapidly becomes more difficult.  One approach to understanding the algebra of covariants is studying its {\it toric degeneration}, see for example \cite{Ca}, \cite{CHV}, \cite{CLO}, \cite{C}, \cite{GL1}, \cite{GL2}, \cite{H2}, \cite{HJLTW}, \cite{Ka}, \cite{KL} and \cite{KM}. To describe the toric degeneration of the algebra $\cal Q$ of covariants, one considers the notion of leading terms with respect to a certain natural ordering on monomials. The leading terms of highest weight polynomials form a semigroup, since the leading term of a product of two polynomials is the product of leading terms of the factors. The algebra generated by the leading terms of  elements of $\cal Q$ is the toric degeneration of $\cal Q$.

In this paper we first revisit, in Section 2, the well studied case of three qubits, i.e., we consider the $G=GL(2)\times GL(2)\times GL(2)$-action on $\bbC^2\otimes\bbC^2\otimes \bbC^2$. Using some ideas of Howe, in particular the use of leading terms, we write very explicitly not only generators and (one) relation defining the covariant algebra $\cal Q$, but also a basis of $\cal Q$. 
As is well known, $\cal Q$ has six generators: one in degree 1, three in degree 2, one in degree 3 and one in degree 4. These generators are denoted respectively by $x_{111},d_{12},d_{13},d_{23},f_3$ and $f_4$, and
they satisfy one relation,
\[
f_3^2=x_{111}^2f_4 -  4d_{12}d_{13}d_{23}.
\]
A basis is then given by monomials in these generators, with $f_3$ appearing only with exponents 0 or 1. 
To prove that these monomials indeed span $\cal Q$, we use some results about adding polynomials over integral points in a polytope, which seem to be interesting in their own right.

We then continue to describe orbits of $G$ in $\bbC^2\otimes\bbC^2\otimes \bbC^2$; this was previously done by Brylinski \cite{B}, but we add some information about equations defining these orbits, which are given in terms of our generators and the irreducible $G$-representations they generate. Along each orbit there is a fixed rank, or degree of entanglement, of its elements. 

In Section 3, we describe in detail the toric degeneration of the algebra $\cal Q$ of covariants attached to three qubits. The semigroup consisting of leading terms of the basis elements is generated by the six leading terms of the generators. We show that this semigroup is affine, as it is a subsemigroup of $\bbZ_+^5$ generated by the six vectors corresponding to the generators. We also show that this semigroup is not saturated, 
hence the corresponding toric variety is not normal.

In Section 4 we study the covariant algebra attached to the $GL(2)\times GL(2)\times GL(3)$-action on the polynomial algebra $\cal P=\cal P(\bbC^2\otimes\bbC^2\otimes \bbC^3)$. This is already much more challenging computationally, and we use computers to obtain all highest weight vectors up to degree 12, and identify the generators and relations. There are thirteen generators, in degrees 1,2,2,2,3,3,4,4,4,5,6,6,6, and there are 20 relations among these generators, all showing up in degrees up to 12. Qualitatively, each of the relations says that a certain product of two generators (or the square of a generator) is redundant, and this enables us to write down a basis of the algebra of covariants. The method of proving that we indeed obtained all generators and relations, and wrote down a (full) basis, is the same as for the case of three qubits, but the arguments are much longer and more complicated. Once we have a basis, it is easy to express the multiplicities of the irreducible components of $\cal P$, as numbers of nonnegative solutions to certain systems of linear equations.

Finally, in Section 5 we explicitly describe the toric degeneration of the algebra of covariants for the ``223 case". The corresponding semigroup $S_{223}$ is generated by the leading terms of the generators, and can be identified with the subsemigroup of $\bbZ_+^7$ generated by the 13 vectors corresponding to the generators. In particular, this semigroup is affine, but as in the 222 case, it is not saturated,
so the toric variety defined by $S_{223}$ is not normal.

We emphasize that in our approach we get to describe all highest weight vectors along with their leading terms very explicitly, which gives more refined information than multiplicity alone.

\section{Covariant algebra attached to the space of three-qubit states}

\subsection{Preliminaries} We consider the action of 
\[
G=G_1\times G_2\times G_3=\GL_2(\bbC)\times \GL_2(\bbC)\times\GL_2(\bbC)
\]
on 
\[
V=V_1\otimes V_2\otimes V_3=\bbC^2\otimes\bbC^2\otimes\bbC^2,
\]
where $V_i=\bbC^2$ is the standard module for $G_i$.
Our main objective is to decompose the algebra  $\Pol=\Pol(V)$ of polynomial functions on $V$ into irreducible $G$-modules.

The vectors
\[
e_i\otimes e_j\otimes e_k,\qquad i,j,k=1,2,
\]
where $e_1,e_2$ denote the standard basis for $\bbC^2$, form a basis for $V$. Let $x_{ijk}$ denote the coordinate functions with respect to this basis. Then $\Pol$ is the algebra of polynomials in the variables $x_{ijk}$. 

If we denote (as usual) by $e_{ij}$ the matrix with the $ij$ entry equal to 1 and other entries equal to 0, then the action of the Lie algebra $\frgl_2(\bbC)$ on $\bbC^2$ is given by
\begin{eqnarray*}
&& e_{11}e_1=e_1,\quad e_{11}e_2=0,\quad e_{22}e_1=0,\quad e_{22}e_2=e_2,\\
&& e_{12}e_1=0,\quad e_{12}e_2=e_1,\quad e_{21}e_1=e_2,\quad e_{21}e_2=0.
\end{eqnarray*}
In particular, $e_1$ is of weight $(1,0)$ for the diagonal matrices $h_1=e_{11}$ and $h_2=e_{22}$, and it is the highest weight vector of the module $\bbC^2$ with respect to the standard choice of positive root vector, $e_{12}$. Furthermore, $e_2$ is of weight $(0,1)$ and it is the lowest weight vector of $\bbC^2$.

This immediately implies that $V$ is an irreducible $G$-module with highest weight $(1,0)\otimes(1,0)\otimes(1,0)=(1,0\bbar 1,0\bbar 1,0)$, with highest weight vector $x_{111}$. The action of the $r$th $\frgl_2$-factor, $\frgl_2^{(r)}$, in the basis $x_{ijk}$ is such that $e_{12}^{(r)}$ changes the $r$th index $2$ to $1$, and annihilates basis elements with the $r$th index equal to 1, while $e_{21}^{(r)}$ changes the $r$th index 1 to 2 and annihilates basis elements with the $r$th index equal to 2. For example, $e_{12}^{(1)}$ sends $x_{121}$ to 0 and $x_{212}$ to $x_{112}$. It is also easy to read off the weights; for example, $x_{221}$ has weight $(0,1\bbar 0,1\bbar 1,0)$. 

We denote by $N_i$ the group of unipotent upper triangular matrices in $G_i$, and set $N=N_1\times N_2\times N_3\subset G$. Clearly, a polynomial $f\in \Pol$ is a highest weight vector if and only if it is invariant for $N$, and a joint eigenvector for the action of the maximal torus $T_1\times T_2\times T_3$ in $G$, where $T_r$ denotes the subgroup of $G_r$ consisting of diagonal matrices.

Passing to the Lie algebra, $f\in \Pol$ is a highest weight vector if and only if it is a weight vector, and annihilated by all $e_{12}^{(r)}$, $r\in\{1,2,3\}$.

For example, $x_{111}\in\Pol$ is a highest weight vector of weight $(1,0\bbar 1,0\bbar 1,0)$, while 
$x_{111}^a$ is a highest weight vector of weight $(a,0\bbar a,0\bbar a,0)$, for any $a\in\bbZ_+$.
Note that the space of $N$-invariants forms a subalgebra of $\Pol$; we call this subalgebra the covariant algebra. Moreover, a product of any two highest weight vectors is again a highest weight vector, and the weight of the product is the sum of the weights of factors.

\subsection{Generators} It is clear that $x_{111}$ is up to factor the only highest weight vector in $\Pol^1=V$; this will be the first generator of our covariant algebra. In degree 2, we follow the idea of Howe's \cite[Theorem 2.1.2]{H} for  the two-factor case, and define
\begin{eqnarray}
\label{dij 222}
&&d_{12}=\left|\begin{matrix} x_{111} & x_{121}\cr x_{211} & x_{221}\end{matrix}\right| = x_{111}x_{221}-x_{121}x_{211};\\
\nonumber
&&d_{13}=\left|\begin{matrix} x_{111} & x_{112}\cr x_{211} & x_{212}\end{matrix}\right| = x_{111}x_{212}-x_{112}x_{211};\\
\nonumber
&&d_{23}=\left|\begin{matrix} x_{111} & x_{112}\cr x_{121} & x_{122}\end{matrix}\right| = x_{111}x_{122}-x_{112}x_{121}.
\end{eqnarray}
Note that $d_{ij}$ is obtained from  $d=\left|\begin{matrix} x_{11} & x_{12}\cr x_{21} & x_{22}\end{matrix}\right|$ on the indices $i$ and $j$, with the third index kept at 1. 

It is clear that all $d_{ij}$ are highest weight vectors, and that they are linearly independent. Moreover, the representations $\pi_{ij}$ they generate have highest weights
\[
\wt(d_{12})=(1,1\bbar 1,1\bbar 2,0);\qquad \wt(d_{13})=(1,1\bbar 2,0\bbar 1,1);\qquad \wt(d_{12})=(2,0\bbar 1,1\bbar 1,1).
\]
Since the irreducible $GL(2)$-representation with highest weight $(\la_1,\la_2)$ has dimension $\la_1-\la_2+1$, each $\pi_{ij}$ is three-dimensional. Another representation in $\Pol^2$ we know is the one generated by $x_{111}^2$, with weight $(2,0\bbar 2,0\bbar 2,0)$ and dimension $27$. Since $\dim\Pol^2=\binom{2+7}{7}=36=27+3\cdot 3$, we see that we have accounted for all the highest weight vectors in $\Pol^2$.

It will be very important for us to keep track of  leading terms of highest weight polynomials. To define the notion of leading term, we order the variables by lexicographical order of the indices, with taking 1 to be higher than 2. In other words,
\[
x_{111}>x_{112}>x_{121}>x_{122}>x_{211}>x_{212}>x_{221}>x_{222}.
\]
We now order monomials by lexicographical order with respect to the above ordering of the variables. For example, the leading terms of the $d_{ij}$ are
\[
\LT(d_{12})=x_{111}x_{221};\qquad \LT(d_{13})=x_{111}x_{212};\qquad \LT(d_{23})=x_{111}x_{122}.
\]
Proceeding to $\Pol^3$, we already know the ``old" highest weight vectors, $x_{111}^3$, $x_{111}d_{12}$, $x_{111}d_{13}$ and $x_{111}d_{23}$. Their respective
weights are
\[
(3,0\bbar 3,0\bbar 3,0);\qquad (2,1\bbar, 2,1\bbar 3,0);\qquad (2,1\bbar 3,0\bbar 2,1);\qquad (3,0\bbar 2,1\bbar 2,1),
\]
so the dimensions are $64$, $16$, $16$ and $16$. 
Since their leading terms,
\[
x_{111}^3;\qquad x_{111}^2x_{221};\qquad x_{111}^2x_{211};\qquad \text{and}\quad x_{111}^2x_{122}
\]
are different, the polynomials are linearly independent. (This is also clear from the fact that they have different weights.) So the dimensions of the corresponding representations add up, to $64+3\cdot 16=112$.

Since $\dim\Pol^3=\binom{3+7}{7}=120$, we see that we are missing 8 dimensions. It is not difficult, even by hand, to find the missing highest weight vector:
\[
f_3= x_{111}^2x_{222}-x_{111}x_{112}x_{221}-x_{111}x_{121}x_{212}-x_{111}x_{122}x_{211}+2x_{112}x_{121}x_{211}.
\]
Its weight and leading term are
\[
\wt(f_3)=(2,1\bbar 2,1\bbar 2,1);\qquad \LT(f_3)=x_{111}^2x_{222}.
\]
Since $\LT(f_3)$ is different from $\LT(x_{111}^3)$ and $\LT(x_{111}d_{ij})$, $f_3$ is linearly independent from $x_{111}^3$ and $x_{111}d_{ij}$, and since the dimension of the representation generated by $f_3$ is 8, we have exhausted $\Pol^3$. 

Passing to $\Pol^4$, we list the ``old" highest weight vectors as $x_{111}^4$, $x_{111}^2d_{ij}$, $d_{ij}d_{kl}$ and $x_{111}f_3$, and patiently compute corresponding weights, leading terms, and dimensions. It turns out that we are missing a one-dimensional representation. One finds the corresponding highest weight polynomial is
\begin{eqnarray*}
& f_4 =& x_{111}^2 x_{222}^2 -2x_{111}x_{112}x_{221}x_{222} -2x_{111}x_{121}x_{212}x_{222}-2x_{111}x_{122}x_{211}x_{222}
+\\
&& 4x_{111}x_{122}x_{212}x_{221}+x_{112}^2 x_{221}^2+4x_{112}x_{121}x_{211}x_{222}-2x_{112}x_{121}x_{212}x_{221}-\\&& 2x_{112}x_{122}x_{211}x_{221}+x_{121}^2x_{212}^2-2x_{121}x_{122}x_{211}x_{212}+x_{122}^2x_{211}^2,
\end{eqnarray*}
with
\[
\wt(f_4)=(2,2\bbar 2,2\bbar 2,2)\qquad\text{and}\qquad \LT(f_4)=x_{111}^2x_{222}^2.
\]
\begin{rem} The less obvious generators $f_3$ and $f_4$ can be obtained as highest weight vectors of the PRV components in $V\otimes\pi_{12}$ respectively $\pi_{12}\otimes\pi_{12}$. (Recall that the PRV (Parthasarathy-Ranga Rao-Varadarajan) component of the tensor product of finite-dimensional representations $V_\la$ and $V_\mu$ of a reductive group with highest weights $\la$ respectively $\mu$ is the unique subrepresentation of $V_\la\ot V_\mu$ with extremal weight equal to the sum of $\la$ and the lowest weight of $V_\mu$; see \cite{PRV}.)
	
Indeed, if we denote 
\begin{eqnarray*}
&&\omega_{12}=e_{21}^{(3)}d_{12}=x_{111}x_{222}+x_{112}x_{221}-x_{121}x_{212}-x_{122}x_{211};\\
&&\eta_{12}=\frac{1}{2}e_{21}^{(3)}\omega_{12}=x_{112}x_{222}-x_{122}x_{212}=
\left|\begin{matrix} x_{112}& x_{122}\cr x_{212}& x_{222}\end{matrix}\right|,
\end{eqnarray*}
then one checks
\[
f_3=x_{111}\omega_{12}-2x_{112}d_{12};\qquad f_4=\omega_{12}^2-4d_{12}\eta_{12},
\]
and these are the highest weight vectors of the PRV components mentioned above.
\end{rem}
As we shall see below, $x_{111},d_{ij},f_3$ and $f_4$ generate the algebra of $N$-invariants in $\Pol$.

\subsection{Relation} 
 Looking at the leading terms, we observe that $f_3^2$ and $x_{111}^2f_4$ have the same leading term, $x_{111}^4x_{222}^2$. This means that $f_3^2-x_{111}^2f_4$ can be expressed as a combination of monomials in the generators with lower leading terms; in fact, one checks that
\eq
\label{rel 222}
f_3^2=x_{111}^2f_4-4d_{12}d_{13}d_{23}.
\eeq
As we shall see below, this is the only relation among our generators, in the sense that all other relations are multiples of this one.

\subsection{Basis}
We can now prove the following well known result (refs...):

\begin{prop}\label{basis 222} 
The monomials
\eq
\label{mono 222}
x_{111}^a d_{12}^b d_{13}^c d_{23}^d f_3^e f_4^g,\qquad a,b,c,d,g\in\bbZ_+,\ e=0,1
\eeq
form a vector space basis for the algebra of $N$-invariants in $\Pol$. In particular, this algebra is generated by $x_{111},d_{12},d_{13},d_{23},f_3$ and $f_4$, with the only relation being \eqref{rel 222}.
\end{prop}

\pf To see that the monomials \eqref{mono 222} are linearly independent, we show that their leading terms are different. Indeed, the leading term of \eqref{mono 222} is 
\eq
\label{LT 222} x_{111}^{a+b+c+d+2e+2g} x_{221}^b x_{212}^c x_{122}^d x_{222}^{e+2g}.
\eeq
This leading term obviously determines $b$, $c$ and $d$, as exponents of $x_{221}$, $x_{212}$ and $x_{122}$ respectively. The exponent of $x_{222}$ determines $e$ and $g$, since $e\in\{0,1\}$. Finally, $a$ is then also determined. 

To see that the monomials \eqref{mono 222} also span the space of $N$-invariants in $\Pol$, we check that for every degree $N$, the sum of dimensions of representations corresponding to the monomials \eqref{mono 222} of degree $N$ equals $\dim\Pol^N=\binom{N+7}{7}$. To see this, we read off
the weight $(\la_1,\la_2\bbar \mu_1,\mu_2\bbar\nu_1,\nu_2)$ of \eqref{mono 222}  from its leading term \eqref{LT 222}, as
\begin{multline}
	\label{h wts}
	(a+b+c+2d+2e+2g,b+c+e+2g\bbar  a+b+2c+d+2e+2g,b+d+e+2g\bbar \\a+2b+c+d+2e+2g,c+d+e+2g).
\end{multline} 
This is obtained by counting the number of 1s and 2s in each of the indices; the number of 1s is the first coordinate of the corresponding component of the weight $(\la\bbar\mu\bbar\nu)$, while the number of 2s is the second coordinate.
We also note that the degree of
\eqref{mono 222}   is
\eq
\label{deg}
N=a+2b+2c+2d+3e+4g=\la_1+\la_2=\mu_1+\mu_2=\nu_1+\nu_2.
\eeq

From the highest weight \eqref{h wts} we see that the dimension of the corresponding representation is
\[
(a+2d+e+1)(a+2c+e+1)(a+2b+e+1).
\]
So it suffices to prove that
\eq
\label{eqn N}
\sum_{a+2b+2c+2d+3e+4g=N} (a+2d+e+1)(a+2c+e+1)(a+2b+e+1) = \binom{N+7}{7}.
\eeq
It is clear that $\binom{N+7}{7}$ is a polynomial in $N$ of degree 7. 
It would be nice if we knew that the left side  of \eqref{eqn N} is also a polynomial in $N$ of degree 7. 
If we knew that, then since two polynomials of degree 7 are equal if their values agree for 8 different values of the argument, it would suffice to check \eqref{eqn N} for $N=0,1,\dots,7$. But we do not know a priori that the left side  of \eqref{eqn N} is a polynomial in $N$ of degree 7; see the next subsection. We are aiming to use Proposition \ref{gen sum} below, but first we eliminate $e$, and write the left side  of \eqref{eqn N} as a sum of two expressions, one for $e=0$ and the other for $e=1$:
\begin{eqnarray*}
&&\sum_{a+2b+2c+2d+4g=N} (a+2d+1)(a+2c+1)(a+2b+1) + \\
&&\sum_{a+2b+2c+2d+4g=N-3} (a+2d+2)(a+2c+2)(a+2b+2).
\end{eqnarray*}
Expressing $a$, we get 
\begin{eqnarray*}
	&&\sum_{2b+2c+2d+4g\leq N} P_1(b,c,d,g,N) + \\
	&&\sum_{2b+2c+2d+4g\leq N-3} P_2(b,c,d,g,N),
\end{eqnarray*}
and this puts us into the setting of Proposition \ref{gen sum}, which tells us is that we should look at congruences of $N$, respectively $N-3$, modulo the least common multiple of the coefficients $2,2,2,4$, which is 4. For each congruence class, each of the above two sums is a polynomial in $N$ of degree at most 7, and this means we should check equality \eqref{eqn N} for 8 numbers in each congruence class modulo 4, that is, for numbers $0,1,\dots,31$. This can be readily done by a computer and this finishes the proof of the proposition modulo Proposition \ref{gen sum}.
\epf

\subsection{Adding polynomials over integral points in certain polytopes} We start by an example showing that it is not clear a priori that the left side of \eqref{eqn N} is a polynomial in $N$.

\begin{ex}
{\rm 
Consider the expression
\eq\label{ex sum 1}
\sum_{\overset{a,b\in\bbZ_+}{a+2b=N}} ab.
\eeq	
We can write $a=N-2b$ and see that the above expression is equal to
\eq\label{ex sum}
\sum_{\overset{b\in\bbZ_+}{2b\leq N}} (N-2b)b.
\eeq
If $N$ is even, $N=2M$, then $2b\leq N$ is equivalent to $b\leq M$, and using the well known formulas $\sum_{b=0}^M b=\frac{M(M+1)}{2}$ and $\sum_{b=0}^M b^2=\frac{M(M+1)(2M+1)}{6}$, we get after some easy manipulations that \eqref{ex sum} is equal to 
\[
\frac{M(M+1)(M-1)}{3}=\frac{N(N+2)(N-2)}{24}.
\]
On the other hand, if $N=2M+1$, then $2b\leq N$ is again equivalent to $b\leq M$, but the relationship between $N$ and $M$ is now different. A similar computation as in the previous case gives
\[
\frac{M(M+1)(2M+1)}{6}=\frac{N(N+1)(N-1)}{24}.
\]
So we see that we got different polynomials for $N$ even respectively odd. On the other hand, if we replace the polynomial $ab$ in \eqref{ex sum 1} by $(a+1)^2$, we again get two polynomials as above, but they both happen to be equal to $\binom{N+3}{3}$.
}
\end{ex}

Motivated by the above example we prove a rather general statement (Proposition \ref{gen sum} below), that will be sufficient for the argument in the proof of Proposition \ref{basis 222}, and also for determining a basis of the covariant algebra for $\bbC^2\ot\bbC^2\ot\bbC^3$ (Theorem \ref{basis 223} below.) We start with the following easy and well known lemma:

\begin{lem}\label{lemma poly} 
Suppose that a function $Q:\bbZ_+\to \bbQ$ is such that the function 
\[
P(N)=Q(N)-Q(N-1),\qquad N\in\bbZ_+,
\] 
is a polynomial of degree $d$. (We set $Q(-1)=0$.) Then $Q$ is a polynomial of degree $d+1$.
\end{lem} 
\pf
Let $P(N)=a_dN^d+a_{d-1}N^{d-1}+\dots + a_0$, with $a_d\neq 0$. Then
\[
Q(N)=\sum_{K=0}^N P(K)=a_d\sum_{K=0}^N K^d  +a_{d-1}\sum_{K=0}^N K^{d-1}+\dots+ a_0\sum_{K=0}^N K^0.
\] 
Since for any $r$, $\sum_{K=0}^N K^r$ is a polynomial in $N$ of degree $r+1$, the claim follows.
\epf

\begin{prop}\label{gen sum}
Let $k_1,\dots,k_s\in\bbN$, and let $P=P(x_1,\dots, x_s,N)$ be a  polynomial with rational coefficients of degree $q$. Let
\[
f(N)=\sum_{\overset{x_1,\dots,x_s\in\bbZ_+}{k_1x_1+\dots+k_s x_s\leq N}}P(x_1,\dots,x_s,N),\qquad N\in\bbZ_+.
\]
Let $L$ be the least common multiple of $k_1,\dots,k_s$  and write $N=LM+R$, with $0\leq R<L$. Then for every fixed $R$, $f(N)=f(LM+R)$ is a polynomial in $N$ (or equivalently in $M$), of degree at most $q+s$.
\end{prop}
\pf
We use induction on the number of variables $s$. If $s=1$, then $L=k_1$ and the inequality $k_1x_1\leq N=k_1M+R$ is equivalent to $x_1\leq M$ (since $0\leq R<k_1$). We need to prove that
\[
\phi(M)=\sum_{x_1\in\bbZ_+,\ x_1\leq M}P(x_1,N)=\sum_{x_1\in\bbZ_+,\ x_1\leq M}P(x_1,k_1M+R)
\]
is a polynomial in $M$ of degree $q+1$; then it will also be a polynomial in $N=k_1M+R$ of the same degree, since $k_1$ and $R$ are constants. 
But by Lemma \ref{lemma poly}, it is enough to prove that 
\[
\phi(M)-\phi(M-1)=P(M,k_1M+R)
\]
is a polynomial in $M$ of degree at most $q$, and this is obvious.

Suppose now that $s>1$ and that we know the claim whenever the number of variables is less than $s$. The idea is to eliminate the variable $x_1$ but we should first make its coefficient equal to 1, by dividing by $k_1$. To do this, we have to rewrite the inequality
\eq\label{poly 0}
k_1x_1+\dots+k_sx_s\leq N=LM+R,\qquad R\text{ fixed},\quad 0\leq R<L
\eeq
involved in the definition of $f(N)$.

Let us temporarily fix an $s$-tuple $x_1,\dots,x_s$ involved in \eqref{poly 0}.  
For $j=2,\dots,s$, let $g_j$ be the greatest common divisor of $k_1$ and $k_j$, and set
$l_j=\frac{k_1}{g_j}$ and  $m_j=\frac{k_j}{g_j}$.  Divide $x_j$ by $l_j$, i.e., set
\[
x_j=l_jy_j + r_j,\qquad 0\leq r_j <l_j.
\]
Since $l_jk_j=k_1m_j$ for all $j$, we can rewrite
 \eqref{poly 0} as
\eq\label{poly 1}
k_1x_1+k_1m_2 y_2+\dots + k_1m_s y_s\leq LM+R-\sum_{j=2}^s k_j r_j.
\eeq
We can rewrite the right side of the inequality \eqref{poly 1} as
\[
L(M-C)+R_1,\qquad 0\leq R_1<L;
\]
here $C$ and $R_1$ are constants which depend only on the $r_j$, but not on the $x_j$. More precisely,
\[
C=\left[\frac{\sum k_jr_j -R}{L}\right];\qquad R_1=LC-\sum k_jr_j+R.
\]
Since the left side of \eqref{poly 1} is divisible by $k_1$, and so is $L(M-C)$, we can replace $R_1$ by $k_1R_2$, where $R_2=\left[\frac{R_1}{k_1}\right]$, and get an equivalent inequality. Dividing by $k_1$ we get that \eqref{poly 1} (or \eqref{poly 0}) is equivalent to 
\eq\label{poly 2}
x_1+m_2 y_2+\dots+m_s y_s\leq L_1(M-C)+R_2, \qquad 0\leq R_2< L_1,
\eeq
where we denoted $\frac{L}{k_1}$ by $L_1$.

We now vary $x_1,\dots,x_s$, but keep the same $r_2,\dots r_s$. Note that the $y_j$ are uniquely determined by the $x_j$ and vice versa, and the equations passing between them are linear. We claim that
\eq\label{poly 3}
\sum_{\overset{x_1,\dots,x_s\in\bbZ_+}{x_1+m_2 y_2+\dots+m_s y_s\leq L_1(M-C)+R_2}} P(x_1,\dots,x_s,N)
\eeq
is a polynomial in $M-C$, or equivalently in $M$, or equivalently in $N$, of degree at most $q+s$. If we prove this, then we are done, since we can add up these polynomials over all choices of the $r_j$ (there are at most $(s-1)k_1$ such choices), and get a polynomial in $M-C$ (or in $M$, or in $N$), of degree at most $q+s$.

Using Lemma \ref{lemma poly}, we see it is enough to prove that
\begin{multline}\label{poly 4}
\sum_{\overset{x_1,\dots,x_s\in\bbZ_+}{x_1+m_2 y_2+\dots+m_s y_s= L_1(M-C)+R_2}} P(x_1,\dots,x_s,N) =\\\sum_{\overset{y_2,\dots,y_s\in\bbZ_+}{m_2 y_2+\dots+m_s y_s\leq L_1(M-C)+R_2}} P(N-m_2y_2-\dots -m_sy_s,x_2,\dots,x_s,N)
\end{multline}
 is a polynomial in $M-C$ (or in $M$, or in $N$), of degree at most $q+s-1$. This will immediately follow from the inductive assumption, if we prove that the least common multiple of $m_2,\dots,m_s$, $\LCM(m_2,\dots,m_s)$, is equal to $L_1$ (recall that $R_2\in[0,L_1)$ does not depend on the $x_j$ or $y_j$, only on the $r_j$ which are fixed).
 
 To prove this last claim, it suffices to show that for any prime number $p$, the exponent of $p$ in $L_1$ is equal to the exponent of $p$ in $\LCM(m_2,\dots,m_s)$. Let $a_j$ be the exponent of $p$ in $k_j$; then the exponent of $p$ in $L$ is $\max(a_1,\dots,a_j)$, and the exponent in $L_1$ is
 \eq\label{LCM 1}
 \max(a_1,\dots,a_j)-a_1=\max(0,a_2-a_1,\dots,a_s-a_1).
 \eeq
On the other hand, the exponent of $p$ in $m_j=\frac{k_j}{\gcd(k_1,k_j)}$ is \eq\label{LCM 2}
b_j=a_j-\min(a_1,a_j)=\left\{ \begin{matrix} 0, & \text{ if } a_j\leq a_1 \cr
	a_j-a_1, & \text{ if } a_j > a_1\end{matrix}
\right.
\eeq
It is now clear that the exponent of $p$ in $\LCM(m_2,\dots,m_s)$, i.e.,  $\max(b_2,\dots,b_s)$, is equal to \eqref{LCM 1}. This finishes the proof.
\epf

\subsection{Multiplicities} Proposition \ref{basis 222} immediately implies

\begin{cor}\label{cor mult 222}
The multiplicity of the highest weight
$(\la_1,\la_2\bbar\mu_1,\mu_2\bbar\nu_1,\nu_2)$
in $\Pol$ is equal to the number of solutions $(a,b,c,d,e,g)$, with $a,b,c,d,g\in\bbZ^+$ and $e\in\{0,1\}$,  of the system of equations 
\begin{eqnarray*}
&&a+b+c+2d+2e+2f=\la_1 \\
&&b+c+e+2f=\la_2 \\
&&a+b+2c+d+2e+2f=\mu_1 \\
&&b+d+e+2f=\mu_2 \\
&&a+2b+c+d+2e+2f=\nu_1 \\
&&c+d+e+2f=\nu_2.
\end{eqnarray*}
\end{cor}
It is possible to solve the above system completely explicitly and obtain the following (known) result:

\begin{prop}
\label{mult 222}
The representation of $G$ with highest weight 
\eq\label{h wt 222}
(\la_1,\la_2\bbar\mu_1,\mu_2\bbar\nu_1,\nu_2)
\eeq
appears in $\Pol$ if and only if:
\begin{enumerate}
\item  for some $N\in\bbZ_+$ (the degree of occurence),
\[
\la_1+\la_2=\mu_1+\mu_2=\nu_1+\nu_2=N;\qquad\text{and}
\]
\item $m\leq M$, where
\[
m=\max(0,\mu_2-\la_2,\nu_2-\la_2)\quad \text{and}\quad M= \min\left(\frac{\mu_2+\nu_2-\la_2-\epsilon}{2},\frac{\la_1-\la_2-\epsilon}{2}\right)
\]
with $\epsilon\in\{0,1\}$ being the parity of 
$\mu_2+\nu_2-\la_2$.
\end{enumerate}
If these conditions are satisfied, then the multiplicity of \eqref{h wt 222} in $\Pol$ is equal to $[M]-m+1$.

The corresponding highest weight vectors can be explicitly described as follows. Choose an integer $\delta\in[m,M]$. Let $\epsilon$ be the remainder of $\mu_2+\nu_2-\la_2$ modulo 2. Set $t=\frac{\mu_2+\nu_2-\la_2-\epsilon}{2}$. The corresponding highest weight vector is then
\[
x_{111}^{-2\delta-\epsilon+\la_1-\la_2} d_{12}^{\delta+\la_2-\nu_2} d_{13}^{\delta+\la_2-\mu_2} d_{23}^\delta f_3^\epsilon f_4^{-\delta+t}.
\]
\end{prop}

\begin{rem}
\label{rmk mult}
To explain better the multiplicities present in the decomposition of $\Pol$, we first note that there is no multiplicity up to degree 5, while in degree 6 the representation with highest weight $(4,2\bbar 4,2\bbar 4,2)$ appears twice; one of the highest weight vectors is $x_{111}^2f_4$, while the other is $d_{12}d_{13}d_{23}$. (Note that these are exactly the vectors involved in the relation \eqref{rel 222}, and also that they have different leading terms, $x_{111}^4x_{222}^2$ respectively $x_{111}^3x_{221}x_{212}x_{122}$.)

We claim that all the multiplicities come from this basic one in degree 6. Indeed, the repetitions of the highest weight come with the exponent $\delta$ of $d_{23}$ varying between $m$ and $[M]$. The smallest $\delta=m$ is such that at least one of the exponents of $d_{12},d_{13},d_{23}$ is 0. Then as $\delta$ grows,  the powers of $d_{12}d_{13}d_{23}$ are appearing, while simultaneously  
the exponent of $x_{111}$ decreases by two and the exponent of $f_4$ decreases by one at each step. In other words, as we are adding powers of $d_{12}d_{13}d_{23}$, we are taking out powers of $x_{111}^2f_4$, until we run out of the factors of the form $x_{111}^2f_4$.
\end{rem}

\subsection{$G$-orbits in $V$ }We now briefly explain how the above results relate to the classification of $G$-orbits in $V=\bbC^2\ot\bbC^2\ot\bbC^2$, or rather in the projective space $\bbP(V)$, which was described by 
Brylinski \cite[Proposition 1.7]{B}. We describe these orbits by equations defined in terms of our generators and other vectors in the representations generated by the generators. We note that vanishing of a single polynomial does not make sense along a $G$-orbit unless this polynomial is $G$-invariant, but simultaneous vanishing of all polynomials in a 
$G$-subrepresentation $\pi$ of $\Pol$ does make sense. Of course, vanishing of all polynomials in $\pi$ is equivalent to the vanishing of a basis of $\pi$.

We have already mentioned the representation $\pi_{12}$ generated by the highest weight polynomial $d_{12}$, with basis $d_{12},\omega_{12},\eta_{12}$.  Analogously we define representations
\[
\pi_{13}=\langle d_{13},\omega_{13},\eta_{13}\rangle;\qquad \pi_{23}=\langle d_{23},\omega_{23},\eta_{23}\rangle.
\]
We will say that $\pi_{12}$ vanishes at a point or on a set, and write $\pi_{12}=0$, to indicate vanishing of all polynomials in the representation $\pi_{12}$, or equivalently, vanishing of $d_{12}$, $\omega_{12}$ and $\eta_{12}$. Analogously, we will consider vanishing of $\pi_{13}$ and $\pi_{23}$.

Brylinski \cite{B} described the 6 $G$-orbits in $V$. The smallest orbit is the orbit of decomposable (or non-entangled) states $u\in\bbP(V)$, i.e., those that can be written as $u=v_1\ot v_2\ot v_3$ for some $v_i\in\bbC^2$. One checks that a state $u$ is decomposable if and only if
all of $\pi_{12}$, $\pi_{13}$ and $\pi_{23}$ vanish on $u$. More generally, the rank of a tensor  $u\in\bbC^2\otimes\bbC^2\otimes\bbC^2$ is the minimal number of decomposable tensors that add up to $u$. The rank can be thought of as a measure of entanglement.

We now recall Brylinski's sets $Y_1,Y_2,Y_3$: $Y_1$ is given as the set of states of the form $u=v\ot w$ with $v\in\bbC^2$ and $w\in\bbC^2\ot\bbC^2$; in other words, 
\eq
\label{Y1}
u=(\alpha e_1+\beta e_2)\ot(\gamma e_1\ot e_1+\delta e_1\ot e_2+\eps e_2\ot e_1+\zeta e_2\ot e_2)
\eeq
for some scalars $\alpha,\beta,\gamma,\delta,\eps,\zeta$. $Y_2$ and $Y_3$ are defined analogously, with the role of the first factor played by the second respectively third factor.

One checks that $Y_1$ is the set of states given by the equations $\pi_{12}=0$ and $\pi_{13}=0$. Analogously, $Y_2$ is the set of states given by $\pi_{12}=\pi_{23}=0$ while $Y_3$ is the set of states given by $\pi_{13}=\pi_{23}=0$.

Recall that the generator $f_3$ of the algebra of highest weight vectors in $\Pol$ can be written as
$f_3=x_{111}\omega_{12}-2x_{112}d_{12}$. It is easy to check that $f_3$ can also be expressed as
\[
f_3=x_{111}\omega_{13}-2x_{121}d_{13}=x_{111}\omega_{23}-2x_{211}d_{23}.
\]
We denote by $\rho_3$ the subrepresentation of $\Pol$ generated by $f_3$. It is straightforward to write down a basis of this representation.
\begin{eqnarray}
\label{basis f3}
&& x_{111}\omega_{12}-2x_{112}d_{12},\quad x_{211}\omega_{12}-2x_{212}d_{12},\\ \nonumber &&x_{121}\omega_{12}-2x_{122}d_{12},\quad 
 x_{221}\omega_{12}-2x_{222}d_{12},\\ 
 \nonumber 
&& 2x_{111}\eta_{12}-x_{112}\omega_{12},\quad
 2x_{211}\eta_{12}-x_{212}\omega_{12},\\ \nonumber && 
 2x_{121}\eta_{12}-x_{122}\omega_{12},\quad
2x_{221}\eta_{12}-x_{222}\omega_{12}. 
\end{eqnarray}
Vanishing of the representation $\rho_3$ on a $G$-orbit in $V$ is equivalent to vanishing of all polynomials \eqref{basis f3}. There are analogous bases of $\rho_3$ in terms of $d_{13},\omega_{13},\eta_{13}$ and of $d_{23},\omega_{23},\eta_{23}$.
Using this, it is straightforward to check
\begin{lem}
\label{f3 vs dij}
(1) If a state $u$ is annihilated by any of the $\pi_{ij}$, then $u$ is annihilated by $\rho_3$;

(2) If a state $u$ is annihilated by $\rho_3$,  then $u$ is annihilated by at least two of the $\pi_{ij}$.
\end{lem}

Finally, let us consider the element $f_4=\omega_{12}^2-4d_{12}\eta_{12}$. The polynomial $f_4$ can also be expressed as
\eq
\label{eqns f4}
f_4=\omega_{13}^2-4d_{13}\eta_{13}=\omega_{23}^2-4d_{23}\eta_{23}.
\eeq
This can be proved by direct computations, but one can also note that $f_4$ spans a unique one-dimensional subrepresentation of $\Pol$ of weight $(2,2)\ot(2,2)\ot(2,2)$, which implies the three expressions, being highest weight vectors, must be proportional to each other. It is then enough to check that their leading terms agree, but it is clear they all have the leading term $x_{111}^2x_{222}^2$.

It is also easy to check that our $f_4$ is the same as the polynomial $D$ of \cite[\S 1.4]{B}.

It is clear from \eqref{eqns f4} that vanishing of any one of the $\pi_{ij}$ implies vanishing of $f_4$. It is then clear from Lemma \ref{f3 vs dij} that vanishing of the representation $\rho_3$ also implies vanishing of $f_4$ (or equivalently vanishing of the representation $\rho_4$ spanned by $f_4$). We are now in a position to compare our results with the $G$-orbit decomposition described in \cite[Proposition 1.7]{B}.

\begin{prop} The six (nonzero) $G$-orbits in $V$, or in $\bbP(V)$, are described as follows:

(1) The open orbit is given by $f_4\neq 0$. It consists of rank two tensors, and it can be represented by the state $e_1\ot e_1\ot e_1+e_2\ot e_2\ot e_2$.

(2) The next orbit, which is open in the zero set of $f_4$, is given by $f_4=0$, $\rho_3\neq 0$. It consists of rank three tensors, and it can be represented by the state $e_1\ot e_1\ot e_1+e_2\ot e_1\ot e_2+e_2\ot e_2\ot e_1$.

(3) The next three orbits are given by $\pi_{12}=\pi_{13}=0$, $\pi_{23}\neq 0$, respectively $\pi_{12}=\pi_{23}=0$, $\pi_{13}\neq 0$ respectively $\pi_{13}=\pi_{23}=0$, $\pi_{12}\neq 0$. Each of these orbits consists of rank two tensors, and they can be represented by the states $e_1\ot e_1\ot e_2+e_1\ot e_2\ot e_1$, respectively $e_1\ot e_1\ot e_2+e_2\ot e_1\ot e_1$, respectively $e_1\ot e_2\ot e_1+e_2\ot e_1\ot e_1$.

(4) The orbit of decomposable (rank one) tensors is given by $\pi_{12}=\pi_{13}=\pi_{23}= 0$ and it can be represented by the state $e_1\ot e_1\ot e_1$.
\end{prop}

The natural next step would be to classify the three-qubit states, i.e., to describe the much finer decomposition of $V$ or $\bbP(V)$ into orbits for the group $U(2)\times U(2)\times U(2)$, or for the group $SU(2)\times SU(2)\times SU(2)$. For this, we refer the reader to \cite[Proposition 1.8]{B} or to \cite{MW1}.

\section{Toric degeneration for the covariant algebra attached to three qubits}

Recall that we have described the generators $x_{111},d_{12},d_{13},d_{23},f_3$ and $f_4$ for the (covariant) algebra $\cA_{222}$ of $N$-invariants in $\Pol=\Pol(\bbC^2\ot\bbC^2\ot\bbC^2)$, along with their leading terms
\eq\label{LTs 222}
x_{111},\ x_{111}x_{221},\ x_{111}x_{212},\ x_{111}x_{122},\ x_{111}^2x_{222},\ x_{111}^2x_{222}^2.
\eeq
It is clear that all the leading terms \eqref{LT 222} of monomials in $\cA_{222}$ define a semigroup (with respect to multiplication). Denote this semigroup by $S_{222}$. It is also clear that $S_{222}$ is generated by the leading terms \eqref{LTs 222}. The corresponding semigroup algebra, also called the initial algebra of $\cA_{222}$, and denoted by $\ini(\cA_{222})$, is therefore finitely generated, and hence is a toric degeneration of $\cA_{222}$, as follows from \cite{CHV} (see also \cite[Proposition 2.5.2]{HJLTW}).

We claim that the semigroup $S_{222}$ is affine, and therefore the toric variety $\Spec_{\ini(\cA_{222})}$ is an affine toric variety. To see this, we first notice that there are five variables, 
$x_{111}$, $x_{122}$, $x_{212}$, $x_{221}$ and $x_{222}$ 
involved in the leading terms \eqref{LTs 222}. Looking at the exponents of these variables, we see that $S_{222}$ can be identified with the subsemigroup of $\bbZ_+^5$ generated by the following 6 vectors:
\[
(1,0,0,0,0),\ (1,0,0,1,0),\ (1,0,1,0,0),\ (1,1,0,0,0),\ (2,0,0,0,1)\ \text{ and } (2,0,0,0,2).
\]
Another way to view the semigroup $S_{222}$ is as the lattice $L$ in $\bbZ_+^5$
generated by $(1,0,0,0,0)$, $(1,0,0,1,0)$, 
$(1,0,1,0,0)$, $(1,1,0,0,0)$  and $ (2,0,0,0,2)$, extended by the element $t=(2,0,0,0,1)$, which is not in $L$ but satisfies
\[
2t=2\cdot (1,0,0,0,0) + (2,0,0,0,2) \in L.
\] 
Note that the semigroup $S_{222}$ is not saturated. Namely,  the group generated by $S_{222}$ contains the element
\[
f=(1,0,0,0,1)=(2,0,0,0,1)-(1,0,0,0,0),
\] 
which is not in $S_{222}$, but $2f$ is in $S_{222}$. It follows that the toric variety defined by $S_{222}$ is not normal; see for example \cite[Theorem 1.3.5]{C}.

\section{Covariant algebra attached to the space of  $(2,2,3)$-qudit states}

\subsection{Computer program} In this section, we study the action of $G=\GL_2(\bbC)\times \GL_2(\bbC)\times\GL_3(\bbC)$ on the  hybrid qudit space
$\bbC^2\otimes\bbC^2\otimes\bbC^3$. 
While the computations in the case of 
three-qubit states can be relatively easily done by hand, they are becoming rather heavy in the present case. Therefore we made a computer program which finds highest weight vectors of a given weight. (We are using the atlas software \cite{atlas}, but all we do could easily be done using better known softwares, like Mathematica.) 

The idea is to write down the matrices of $e_{12}^1$, $e_{12}^2$, $e_{12}^3$ and $e_{23}^3$ in the basis of monomials, put them of top of each other, and then find the kernel. (This kernel will be the intersection of the four kernels.) 

Another thing that the program does is computing the ``old" highest weight vectors, i.e., those that are obtained as products of lower degree ones. This amounts to writing a given weight as a sum of weights appearing in lower degrees, and keeping track of the corresponding highest weight vectors. 

Then we go over degrees $1,2,3,\dots$ and identify the generators, as highest weight vectors that are not ``old".

Finally, to find relations, we look for weights that are obtained as sums of lower degree weights in more than one way. We identify the corresponding monomials (in the generators), and find their leading terms. If two monomials have the same leading term, they are equal up to lower monomials (i.e., those with smaller leading terms). In the following, we list thus obtained generators and relations. After that we indicate how to obtain a (monomial) basis of the space of highest weight vectors.

\subsection{Generators}\label{subsec:gens}

\bigskip
\bigskip

\begin{center}
	\begin{tabular}{|c|c|c|c|c|}
		\hline
		no & deg & name & leading term & weight  \\ \hline\hline
		0 & 1 & $x_{111}$ & $x_{111}$ & $(10 \sbar 10 \sbar 100)$ \\ \hline
		\hline
		1 & 2 & $d_{12}$ & $x_{111}x_{221}$ & $(11 \sbar 11 \sbar 200)$ \\ \hline
		2 & 2 & $d_{13}$ & $x_{111}x_{212}$ & $(11 \sbar 20 \sbar 110)$ \\ \hline
		3 & 2 & $d_{23}$ & $x_{111}x_{122}$ & $(20 \sbar 11 \sbar 110)$ \\ \hline\hline
		4 & 3 & $f_3$ & $x_{111}^2x_{222}$ & $(21 \sbar 21 \sbar 210)$ \\ \hline
		5 & 3 & $u$ & $x_{111}x_{122}x_{213}$ & $(21 \sbar 21 \sbar 111)$ \\ \hline\hline
		6 & 4 & $f_4$ & $x_{111}^2x_{222}^2$ & $(22 \sbar 22 \sbar 220)$ \\ \hline
		7 & 4 & $v$ & $x_{111}^2x_{212}x_{223}$ & $(22 \sbar 31 \sbar 211)$ \\ \hline
		8 & 4 & $v^*$ & $x_{111}^2x_{122}x_{223}$ & $(31 \sbar 22 \sbar 211)$ \\ \hline\hline
		9 & 5 & $f_5$ & $x_{111}^2x_{122}x_{212}x_{223}$ & $(32 \sbar 32 \sbar 221)$ \\ \hline\hline
		10 & 6 & $\varphi$ & $x_{111}^2 x_{122}x_{212}x_{223}^2$ & $(33 \sbar 33 \sbar 222)$ \\ \hline
		11 & 6 & $\psi$ & $x_{111}^3 x_{212}x_{222}x_{223}$ & $(33 \sbar 42 \sbar 321)$ \\ \hline
		12 & 6 & $\psi^*$ & $x_{111}^3 x_{122}x_{222}x_{223}$ & $(42 \sbar 33 \sbar 321)$ \\ \hline
	\end{tabular}
\end{center}
\bigskip
\bigskip

The full expressions for some of these generators are as follows:
\begin{eqnarray*}
	& d_{12}=&\left|\begin{matrix} x_{111} & x_{121}\cr x_{211} & x_{221}\end{matrix}\right| \\
	& d_{13}=&\left|\begin{matrix} x_{111} & x_{112}\cr x_{211} & x_{212}\end{matrix}\right| \\
	& d_{23}=&\left|\begin{matrix} x_{111} & x_{112}\cr x_{121} & x_{122}\end{matrix}\right| \\
	& f_3 =& x_{111}^2x_{222}-x_{111}x_{112}x_{221}-x_{111}x_{121}x_{212}-x_{111}x_{122}x_{211} +2 x_{112}x_{121}x_{211} \\
	& u =& \left|\begin{matrix} x_{111} & x_{112} & x_{113} \cr
		x_{121} & x_{122} & x_{123}\cr
		x_{211} & x_{212} & x_{213}
	\end{matrix}\right| 
	\\
	& f_4 =& x_{111}^2 x_{222}^2 -2x_{111}x_{112}x_{221}x_{222} -2x_{111}x_{121}x_{212}x_{222}-2x_{111}x_{122}x_{211}x_{222}
	+\\
	&& 4x_{111}x_{122}x_{212}x_{221}+x_{112}^2 x_{221}^2+4x_{112}x_{121}x_{211}x_{222}-2x_{112}x_{121}x_{212}x_{221}-\\&& 2x_{112}x_{122}x_{211}x_{221}+x_{121}^2x_{212}^2-2x_{121}x_{122}x_{211}x_{212}+x_{122}^2x_{211}^2\\
	& v =& \left|\begin{matrix} 0&x_{111} & x_{112} & x_{113} \cr
		0&x_{211} & x_{212} & x_{213}\cr
		x_{111}&x_{121} & x_{122} & x_{123}\cr
		x_{211}&x_{221} & x_{222} & x_{223}
	\end{matrix}\right| \\
	& v^* =& \left|\begin{matrix} 0&x_{111} & x_{112} & x_{113} \cr
		0&x_{121} & x_{122} & x_{123}\cr
		x_{111}&x_{211} & x_{212} & x_{213}\cr
		x_{121}&x_{221} & x_{222} & x_{223}
	\end{matrix}\right|
\end{eqnarray*}

Note how $v^*$ is obtained from $v$ by switching the first two indices, as one can predict from the highest weights.

We do not write $f_5$ (29 terms), $\varphi$ (66 terms), $\psi$ or $\psi^*$ (44 terms each). These are very big expressions, and even if they do have some nice determinant descriptions, we do not think we can see that by staring at the expressions. (We were able to do the staring for $v$ and $v^*$, encouraged by the fact that they have 12 terms, which is half of $4!$.)

\subsection{Relations}\label{subsec:rels}
We do not write the full relations, but only indicate which monomials in the generators are equal modulo lower terms (i.e., terms with smaller leading terms). So instead of $=$, we write $\equiv$. This suffices for the purpose of writing a basis of the algebra of $N$-invariants.

\bigskip

\begin{center}
	\begin{tabular}{|c|c|c|}
		\hline
		no & deg & relation  \\ \hline\hline
		R1 & 6 & $f_3^2\equiv x_{111}^2f_4$ \\ \hline
		\hline
		R2 & 6 & $d_{23}v\equiv x_{111}f_5$  \\ \hline
		R3 & 6 & $d_{13}v^*\equiv x_{111}f_5$  \\ \hline\hline
		R4 & 7 & $f_3 v\equiv x_{111}\psi$  \\ \hline
		R5 & 7 & $f_3 v^*\equiv x_{111}\psi^*$  \\ \hline\hline
		R6 & 8 & $f_3 f_5\equiv d_{23}\psi$  \\ \hline
		R7 & 8 & $d_{13}\psi^*\equiv d_{23}\psi$  \\ \hline
		R8 & 8 & $vv^*\equiv x_{111}^2\varphi$ \\ \hline\hline
		R9 & 9 & $v f_5\equiv x_{111}d_{13}\varphi$ \\ \hline
		R10 & 9 & $v^* f_5\equiv x_{111}d_{23}\varphi$  \\ \hline
		R11 & 9 & $f_3\psi\equiv x_{111}f_4 v$  \\ \hline
		R12 & 9 & $f_3 \psi^*\equiv x_{111}f_4 v^*$  \\ \hline\hline
		R13 & 10 & $f_5^2 \equiv d_{13}d_{23}\varphi$  \\ \hline
		R14 & 10 & $v^*\psi\equiv x_{111}f_3\varphi$  \\ \hline   
		R15 & 10 & $v\psi^*\equiv x_{111}f_3\varphi$  \\ \hline \hline
		R16 & 11 & $f_5\psi\equiv d_{13}f_3\varphi$  \\ \hline 
		R17 & 11 & $f_5\psi^*\equiv d_{23}f_3\varphi$  \\ \hline \hline
		R18 & 12 & $\psi^2\equiv f_4v^2$  \\ \hline 
		R19 & 12 & $(\psi^*)^2\equiv f_4(v^*)^2$  \\ \hline 
		R20 & 12 & $\psi\psi^*\equiv x_{111}^2f_4\varphi$  \\ \hline 
	\end{tabular}
\end{center}


\subsection{Proof that we have a full description}

We now use the relations to eliminate some of the monomials. For the remaining monomials, we check that they all have different leading terms, and therefore are linearly independent. Finally, for each degree $N$ we add up the dimensions of modules corresponding to the remaining monomials of degree $N$, and  see that we obtain 
\[
\binom{N+11}{11}=\dim\cal P^N(\bbC^2\otimes\bbC^2\otimes\bbC^3).
\]

The scheme is as follows:

By the relations (R1), (R6) and (R13), $f_3$ and $f_5$ appear with exponents 0 or 1, but not both with exponent 1. If $f_3$ appears and $f_5$ does not, then using the relations we see that $v,v^*,\psi$ and $\psi^*$ can not appear, and if they do not appear then all the relations become irrelevant. Thus the surviving monomials containing $f_3$ are
\eq
\label{mono 1}
x_{111}^ad_{12}^bd_{13}^cd_{23}^df_3u^ff_4^g\varphi^k.
\eeq
The corresponding leading term is
\eq\label{LT 1}
x_{111}^{a+b+c+d+2+f+2g+2k}x_{122}^{d+f+k}x_{212}^{c+k}x_{213}^{f}x_{221}^bx_{222}^{1+2g}x_{223}^{2k}.
\eeq
We see that going from right to left we can reconstruct all the exponents in \eqref{mono 1}, so all monomials \eqref{mono 1} have different leading terms.

If $f_3$ does not appear and $f_5$ does, the relations again imply that
$v,v^*,\psi$ and $\psi^*$ do not appear, and if they do not appear then all the relations become irrelevant. Thus we are led to monomials
\eq\label{mono 2}
x_{111}^ad_{12}^bd_{13}^cd_{23}^du^ff_4^gf_5\varphi^k.
\eeq
The corresponding leading term is
\eq\label{LT 2}
x_{111}^{a+b+c+d+f+2g+2+2k}x_{122}^{d+f+1+k}x_{212}^{c+1+k}x_{213}^{f}x_{221}^bx_{222}^{2g}x_{223}^{1+2k}.
\eeq
We again see that going from right to left we can reconstruct all the exponents in \eqref{mono 2}, so all monomials \eqref{mono 2} have different leading terms. Moreover, none of the leading terms in \eqref{LT 2} can be equal to any of the leading terms in \eqref{LT 1}, since they have opposite parity of the exponent of $x_{223}$ (and also of $x_{222}$).

Let us now assume that $f_3$ and $f_5$ do not appear. This makes relations (R1), (R4), (R5), (R6), (R9), (R10), (R11), (R12), (R13), (R16) and (R17) irrelevant. 
We make a discussion with respect to $v$ and $v^*$; they can not both appear by (R8), so we have three cases: $v$ appears and $v^*$ does not, $v^*$ appears and $v$ does not, or none of them appears. 

Suppose that $v$ appears and $v^*$ does not. Then $d_{23}$ does not appear by (R2), and $\psi^*$ does not appear by (R15). With these conditions all relations become irrelevant, except for (R18) which says that $\psi$ appears with exponent at most 1. Thus we are led to monomials of the form
\eq\label{mono 3}
x_{111}^ad_{12}^bd_{13}^c u^f f_4^g v^h \varphi^k\psi^l,\qquad h\geq 1,\ l=0\text{ or }1.
\eeq
The corresponding leading term is
\eq\label{LT 3}
x_{111}^{a+b+c+f+2g+2h+2k+3l}x_{122}^{f+k}x_{212}^{c+h+k+l}x_{213}^{f}x_{221}^bx_{222}^{2g+l}x_{223}^{h+2k+l}.
\eeq
We see immediately that the leading term determines $f$ and $b$, and also $g$ and $l$ since $l$ is 0 or 1. Then the exponent of $x_{122}$ determines $k$, and the exponent of $x_{223}$ determines $h$. Now the exponent of $x_{212}$ determines $c$ and finally the exponent of $x_{111}$ determines $a$. So all the leading terms \eqref{LT 3} are different. 

Suppose that \eqref{LT 3} is equal to \eqref{LT 1} with $a,\dots,k$ replaced by $a',\dots,k'$. Then we see that $f'=f$, so the exponent of $x_{122}$ gives $k=d'+k'$, hence $k\geq k'$. On the other hand, since $h\geq 1$, the exponent of $x_{223}$ gives $k'>k$, a contradiction. We similarly see that \eqref{LT 3} can not be equal to a leading term of the form \eqref{LT 2}.

If a monomial does not contain $f_3,f_5,v$ but does contain $v^*$, then the relations imply that it is of the form
\eq\label{mono 4}
x_{111}^ad_{12}^bd_{23}^d u^f f_4^g {v^*}^i \varphi^k{\psi^*}^m,\qquad i\geq 1,\ m=0\text{ or }1.
\eeq
The corresponding leading term is
\eq\label{LT 4}
x_{111}^{a+b+d+f+2g+2i+2k+3m}x_{122}^{d+f+i+k+m}x_{212}^kx_{213}^{f}x_{221}^bx_{222}^{2g+m}x_{223}^{i+2k+m}.
\eeq
Similar considerations as above show that the leading terms \eqref{LT 4} are different from each other, and also from the leading terms of the form \eqref{LT 1}, \eqref{LT 2} or \eqref{LT 3}.

We now assume that a monomial does not contain $f_3,f_5,v,v^*$, but does contain $d_{13}$. The relations then imply it is of the form
\eq\label{mono 5}
x_{111}^ad_{12}^bd_{13}^cd_{23}^d u^f f_4^g  \varphi^k\psi^l,\qquad c\geq 1,\ l=0\text{ or }1.
\eeq
The corresponding leading term is
\eq\label{LT 5}
x_{111}^{a+b+c+d+f+2g+2k+3l}x_{122}^{d+f+k}x_{212}^{c+k+l}x_{213}^{f}x_{221}^bx_{222}^{2g+l}x_{223}^{2k+l}.
\eeq
As before, we see that these leading terms are different from each other and from the previously considered ones.

We next assume that a monomial does not contain $f_3,f_5,v,v^*,d_{13}$, but does contain $\psi$. The relations then imply it is of the form
\eq\label{mono 6}
x_{111}^ad_{12}^bd_{23}^d u^f f_4^g  \varphi^k\psi.
\eeq
The corresponding leading term is
\eq\label{LT 6}
x_{111}^{a+b+d+f+2g+2k+3}x_{122}^{d+f+k}x_{212}^{k+1}x_{213}^{f}x_{221}^bx_{222}^{2g+1}x_{223}^{2k+1}.
\eeq
As before, we see that these leading terms are different from each other and from the previously considered ones.

Finally, let us assume that a monomial does not contain $f_3,f_5,v,v^*,d_{13}$ or $\psi$. The relations then imply it is of the form
\eq\label{mono 7}
x_{111}^ad_{12}^bd_{23}^d u^f f_4^g  \varphi^k{\psi^*}^m,\qquad m=0\text{ or }1.
\eeq
The corresponding leading term is
\eq\label{LT 7}
x_{111}^{a+b+d+f+2g+2k+3m}x_{122}^{d+f+k+m}x_{212}^{k}x_{213}^{f}x_{221}^bx_{222}^{2g+m}x_{223}^{2k+m}.
\eeq
As before, we see that these leading terms are different from each other and from the previously considered ones.

To finish the proof, we let the computer compute the total dimension in degree $N$. It has to make a list of all monomials as above of degree $N$, and then add up the corresponding  dimensions, which are computed by the Weyl dimension formula. This gives the number $\operatorname{totaldim}(N)$, and then we compare it with $\binom{N+11}{11}$. 

For example, \eqref{mono 1} leads to
\[
\sum_{\overset{a,b,c,d,f,g,k\in\bbZ_+}{a+2b+2c+2d+3f+4g+6k=N-3}}\WDF(\la\bbar\mu\bbar\nu),
\]
where $(\la\bbar\mu\bbar\nu)$ is the weight of the monomial \eqref{mono 1}, 
and $\WDF$ is the Weyl dimension formula polynomial, given by
\begin{multline*}
\WDF(\la_1,\la_2\bbar\mu_1,\mu_2\bbar\nu_1,\nu_2,\nu_3)=\\ \frac{(\la_1-\la_2+1)(\mu_1-\mu_2+1)(\nu_1-\nu_2+1)(\nu_2-\nu_3+1)(\nu_1-\nu_3+2)}{2}
\end{multline*}
As before, we first eliminate $a$ and get
\[
\sum_{\overset{b,c,d,f,g,k\in\bbZ_+}{2b+2c+2d+3f+4g+6k\leq N-3}}\WDF(\la\bbar\mu\bbar\nu),
\]
To use Proposition \ref{gen sum}, we notice that the least common multiple of the coefficients $2,2,2,3,4,6$ is 12, so for each congruence class modulo 12 we get a polynomial of degree 11 (the degree of $\WDF$ is 5, and there are 6 variables, $b,c,d,f,g$ and $k$.)

The situation is analogous for all basic monomials listed above; when we have options for $\psi$ or $\psi^*$ appearing or not, we split the monomials into two cases; sometimes there are only five variables, so the degree is only 10; but when we add up all the polynomials the degree will be 11. 

We conclude that for each congruence class modulo 12 we get a polynomial of degree 11. This means that we have to check the equality
\[
\operatorname{totaldim}(N)=\binom{N+11}{11}
\]
for $12\times 12=144$ values of $N$, i.e., for $N=0,1,\dots,143$. This can easily be done by a computer, and in this way we have proved

\begin{thm} 
	\label{basis 223}
	The algebra $\cA_{223}$ of $N$-invariants in $\Pol$ is generated by the generators described in Subsection \ref{subsec:gens}, subject to relations described in Subsection \ref{subsec:rels}. Monomials in the generators listed in \eqref{mono 1}, \eqref{mono 2}, \eqref{mono 3}, \eqref{mono 4}, \eqref{mono 5}, \eqref{mono 6} and \eqref{mono 7} form a basis of  $\cA_{223}$.	
\end{thm}

\subsection{Multiplicities}
Theorem \ref{basis 223} immediately implies the following corollary:

\begin{cor} Let $\wt_1,\dots,\wt_7$ denote the weights of the monomials described in \eqref{mono 1}, \eqref{mono 2}, \eqref{mono 3}, \eqref{mono 4}, \eqref{mono 5}, \eqref{mono 6} and \eqref{mono 7}, as (linear) functions of the exponents of the generators. For a weight $\la$, and for $i=1,\dots,7$, let $\mult_i$ denote the number of non-negative solutions to the system of equations
\[
\wt_i=\la.
\] 
Then the multiplicity of $\la$ in $\cA_{223}$ equals $\mult_1+\dots+\mult_7$.
\end{cor}

\section{Toric degeneration of the covariant algebra attached to the  hybrid qudit space
$\bbC^2\otimes\bbC^2\otimes\bbC^3$}

Recall that we have described the generators $x_{111},d_{12},d_{13},d_{23},f_3,u,f_4,v,v^*,f_5,\varphi,\psi$ and $\psi^*$ for the (covariant) algebra $\cA_{223}$ of $N$-invariants in $\Pol=\Pol(\bbC^2\otimes\bbC^2\ot\bbC^3)$, along with their leading terms
\begin{eqnarray}\label{LTs 223}
&&x_{111},\\ \nonumber
&&x_{111}x_{221},\\ \nonumber
&&x_{111}x_{212},\\ \nonumber
&&x_{111}x_{122},\\\nonumber
&&x_{111}^2x_{222},\\ \nonumber
&&x_{111}x_{122}x_{213}, \\\nonumber
&&x_{111}^2x_{222}^2,\\ \nonumber
&&x_{111}^2x_{212}x_{223}, \\\nonumber
&&x_{111}^2x_{122}x_{223}, \\ \nonumber
&&x_{111}^2x_{122}x_{212}x_{223}, \\ \nonumber
&&x_{111}^2 x_{122}x_{212}x_{223}^2,\\ \nonumber
&&x_{111}^3 x_{212}x_{222}x_{223}, \\ \nonumber
&&x_{111}^3 x_{122}x_{222}x_{223}.
\end{eqnarray}
It is clear that all the leading terms of monomials in $\cA_{223}$ define a semigroup (with respect to multiplication). Denote this semigroup by $S_{223}$. It is also clear that $S_{223}$ is generated by the leading terms \eqref{LTs 223}. The corresponding semigroup algebra, the initial algebra $\ini(\cA_{223})$ of $\cA_{223}$, is therefore finitely generated, and hence is a toric degeneration of $\cA_{223}$, as follows from \cite{CHV} (see also \cite[Proposition 2.5.2]{HJLTW}). 

We claim that the semigroup $S_{223}$ is affine, and therefore the toric variety $\Spec_{\ini(\cA_{223})}$ is an affine toric variety. To see this, we first notice that there are seven variables, $x_{111},x_{122},x_{212},x_{213},x_{221},x_{222}$ and $x_{223}$, involved in the leading terms \eqref{LTs 223}. Looking at the exponents of these variables, we see that $S_{223}$ can be identified with the subsemigroup of $\bbZ_+^7$ generated by the following 13 vectors:
\begin{eqnarray*}
&&(1,0,0,0,0,0,0),\ (1,0,0,0,1,0,0),\ (1,0,1,0,0,0,0),\ (1,1,0,0,0,0,0),\ (2,0,0,0,0,1,0),\\  &&(1,1,0,1,0,0,0),\ (2,0,0,0,0,2,0),\ (2,0,1,0,0,0,1),\ (2,1,0,0,0,0,1),\\
&&(2,1,1,0,0,0,1),\ (2,1,1,0,0,0,2),\ (3,0,1,0,0,1,1),\ (3,1,0,0,0,1,1).
\end{eqnarray*}
As in the 222 case, the semigroup $S_{223}$ is not saturated. Namely, as before, the group generated by $S_{223}$ contains the element
\[
f=(1,0,0,0,0,1,0)=(2,0,0,0,0,1,0)-(1,0,0,0,0,0,0),
\] 
which is not in $S_{223}$, but $2f$ is in $S_{223}$. It now follows from \cite[Theorem 1.3.5]{C} that the toric variety defined by $S_{223}$ is not normal.


\begin{thebibliography}{HJLTW}

\bibitem[atlas]{atlas} Atlas of Lie Groups and Representations, version 1.2, axis language version 2.0. See www.liegroups.org for more about the software.

\bibitem[BV]{BV} 
V.~Baldoni, M.~Vergne, {\it Multiplicity of compact group representations
and applications to Kronecker coefficients}, arXiv:1506.02472.

\bibitem[BLT]{BLT} E.~Briand, J.-G.~Luque, J.-Y.~Thibon, {\it A complete set of covariants of the four qubit system}

\bibitem[Bri]{Bri} M.~Brion, {\it Introduction to actions of algebraic groups}, Les cours du C.I.R.M., 
Vol. 1 (2010), no. 1, 1--22.

\bibitem[Bry]{B} J-L.~Brylinski, {\it Algebraic measures of entanglement}, Comput. Math. Ser. Chapman \& Hall/CRC, Boca Raton, FL, 2002, 3--23.

\bibitem[BB]{BB} J-L.~Brylinski, R.~Brylinski, {\it Polynomial invariants for qubits}, Comput. Math. Ser. Chapman \& Hall/CRC, Boca Raton, FL, 2002, 277--286.

\bibitem[Ca]{Ca} P.~Caldero, {\it Toric degenerations of	Schubert varieties}, Transform. Groups \textbf{7} (2002), No. 1, 51--60.

\bibitem[CHV]{CHV} A. Conca, J. Herzog and G. Valla,
\textit{SAGBI bases with applications to blow-up algebras}, J. Reine Angew. Math. \textbf{474} (1996), 113--138. 

\bibitem[CLO]{CLO} D.~Cox, J.~Little and D.~O'Shea,
\textit{Ideals, Varieties, and Algorithms. An Introduction to Computational Algebraic Geometry and Commutative Algebra}, Second edition, Undergraduate Texts in Mathematics, Springer-Verlag, New York, 1997.

\bibitem[CLS]{C} D.~Cox, J.~Little and H.~Schenck, {\it Toric varieties}, Graduate Studies in Mathematics 124,  Amer. Math. Soc., Providence, RI, 2011.

  \bibitem[FH]{FH} J.~Fei, C.~Xue, {\it A polyhedral formula for $n\times 2\times 2$ Kronecker coefficients via cluster algebras}, preprint, arXiv:2607.25201.

\bibitem[GL1]{GL1} N. Gonciulea and V. Lakshmibai,
{\it Degenerations of flag and Schubert varieties to toric varieties}, Transform. Groups 1 (1996), No. 3, 215 -- 248.

\bibitem[GL2]{GL2} N.~Gonciulea and V.~Lakshmibai,
\textit{Schubert varieties, toric varieties and ladder determinant	varieties}, Ann. Inst. Fourier (Grenoble) 47 (1997), 1013 -- 1064.

\bibitem[H1]{H} R.~Howe, \textit{Perspectives on invariant theory}, The Schur Lectures, I. Piatetski-Shapiro and S. Gelbart (eds.), Israel Mathematical Conference Proceedings, 1995, 1--182.

\bibitem[H2]{H2} R.~Howe, {\it Weyl chambers and standard monomial theory for poset lattice cones}, Pure Appl. Math. Q. \textbf{1} (2005), no. 1, 227--239.

\bibitem[HJLTW]{HJLTW} R.~Howe, S.~Jackson, S.T.~Lee, E.-C.~Tan, J.~Willenbring, {\it Toric degeneration of branching algebras}, Adv. Math. \textbf{220} (6) (2009), 1809--1841.

\bibitem[Ka]{Ka} K.~Kaveh, {\it SAGBI bases and degeneration of spherical varieties to toric varieties}, Michigan Math. J. \textbf{53} (2005), 
109--121.

\bibitem[KL]{KL} S.~Kim, S.T.~Lee,
\textit{Toric degeneration of algebras of invariants}, Adv. Math. \textbf{459} (2024), 110017.

\bibitem[KM]{KM} M.~Kogan and E.~Miller, {\it
	Toric degeneration of Schubert varieties and Gelfand-Tsetlin polytopes}, Adv. Math. \textbf{193} (2005), 1--17.

\bibitem[La]{La} J.~M.~Landsberg,  
{\it Quantum computation and quantum information: a mathematical perspective}, Graduate Studies in Mathematics \textbf{243}, 2024.

\bibitem[LP]{LP} C. Le Paige, {\it Sur les formes trilin\'eaires}, C. R. Acad. Sci. Paris \textbf{92} (1881), 1103. 

\bibitem[Lu]{Lu} J.-G.~Luque, {\it Invariants des hypermatrices}, habilitation thesis, 2007.

\bibitem[MW1]{MW1} D.~A.~Meyer, N.~Wallach, {\it Invariants for multiple qubits: the case of 3 cubits}, Comput. Math. Ser. Chapman \& Hall/CRC, Boca Raton, FL, 2002, 77--97.

\bibitem[MW2]{MW2} D.~A.~Meyer, N.~Wallach, {\it Global entanglement in multiparticle systems}, J. Math. Phys. \textbf{43} (2002), no. 9, 4273--4278.

\bibitem[PRV]{PRV} K.R.~Parthasarathy, R.~Ranga Rao, V.S.~Varadarajan, \emph{Representations of complex semisimple Lie groups and Lie algebras}, Ann. of Math. \textbf{85} (1967), 383--429.

\bibitem[W]{W} M.~Walter, \textit{Multipartite quantum states
and their marginals}, dissertation, ETH Zurich, 2014, arXiv:1410.6820

\bibitem[WDGC]{WDGC} M.~Walter, B.~Doran,  D.~Gross, M.~Christandl, {\it Entanglement polytopes: multiparticle entanglement from single-particle information}, arXiv:1208.0365

\end{thebibliography}
\end{document}